\documentclass[10pt]{article}
\usepackage{amsmath}
\usepackage{amsfonts}
\usepackage{latexsym}
\usepackage{amssymb}
\usepackage[dvips]{graphics}

\newcommand{\Z}{{\mathbf Z}}

\newcommand{\R}{{\mathbf R}}
\newcommand{\C}{{\mathbf C}}

\newcommand{\cat}{{\rm {cat }}}
\newcommand{\Cat}{\rm {Cat}}
\newcommand{\id}{{\rm {id }}}

\newcommand{\Hom}{{\rm Hom}}
\newcommand{\im}{\rm im}

\newtheorem{theorem}{Theorem}[section]

\newtheorem{lemma}[theorem]{Lemma}

\newtheorem{definition}[theorem]{Definition}

\newtheorem{remark}[theorem]{Remark}

\numberwithin{equation}{section}

\newenvironment{proof}{\par\noindent{\bf Proof.} }{\par\noindent}

\title{Lusternik -- Schnirelman Theory and Dynamics}
\author{Michael Farber\footnote{Partially supported by the US - 
Israel Binational
Science Foundation}}
\date{January 24, 2002}

\begin{document}

\sloppy

\maketitle

\begin{abstract}
In this paper we study a new 
topological invariant $\Cat(X,\xi)$, where $X$ is a finite polyhedron and $\xi\in H^1(X;\R)$ is a real cohomology class.
$\Cat(X,\xi)$ is defined using open covers of $X$ with certain geometric properties; it is 
a generalization of the classical Lusternik -- Schnirelman category. 
We show that $\Cat(X,\xi)$ depends only on the homotopy type of $(X,\xi)$. We prove that $\Cat(X,\xi)$ 
allows to establish a relation between the number of equilibrium states of dynamical systems and their 
global dynamical properties (such as existence of homoclinic cycles and the structure of the set of chain recurrent points). 
In the paper we give  a cohomological lower bound for $\Cat(X,\xi)$, 
which uses cup-products of cohomology classes 
of flat line bundles with monodromy described by complex numbers, which are not Dirichlet units.
\end{abstract}

\section {Introduction}

As is well-known, the geometry of dynamical systems plays an essential role in the classical Lusternik -- Schnirelman and Morse theories.
 Any function determines a flow (the gradient flow), which serves as an important intermediate object connecting the topology of the manifold 
with information about the critical points of the function. Roughly, the gradient flow leads to 
 a partition of the manifold into pieces labeled by the 
critical points.

Gradient systems are very special, and the theory of Conley \cite{C} extends  the methods of the Morse theory 
for much more general dynamical systems. 

At present, the methods of homotopy theory are widely used in dynamics; we refer to \cite{Co} and \cite{RO} for two quite different effective approaches. 

This paper represents an attempt to extend the Lusternik -- Schnirelman theory to flows of a more general nature. 
Our main assumption about a flow on a manifold $M$ is the existence of a Lyapunov closed 1-form $\omega$. This notion 
(which we define in the paper) is a multivalued
generalization of the well-known notion of a Lyapunov function (see \S \ref{final} below). 
We also require that the equilibrium states are isolated in the set of chain recurrent points;
this condition forbids, for example, the flow to have homoclinic orbits and homoclinic cycles. Under these assumptions we show that the number of equilibrium states of the flow
(i.e. the number of zeros of vector field $v$) is at least $\Cat(M, \xi)$. Here $\Cat(M,\xi)$ is a number which depends only
on the homotopy type of $M$
and on the cohomology class $\xi\in H^1(M;\R)$ of the Lyapunov closed 1-form $\omega$
for $v$. For $\xi=0$ the invariant $\Cat(M,\xi)$ coincides with the 
Lusternik -- Schnirelman category $\cat(M)$. 
Note that in the case $\xi=0$ the chain recurrent set coincides with the set of zeros and so our assumptions are
automatically satisfied.

For any closed smooth manifold $M$ and for any nonzero cohomology class $\xi\in H^1(M;\Z)$, $\xi\neq 0$,
there exists a closed 1-form $\omega$ in class $\xi$ having at most one zero 
-- see \cite{F1} for a proof. This Theorem shows that our statements, estimating from below the number of zeros of vector fields 
with Lyapunov closed 1-forms, are incorrect without extra assumptions 
(such as the requirement that the zeros are isolated in the set of all chain recurrent points).
In \cite{F1} these extra requirements were weaker: we only asked about the absence of homoclinic cycles. On the other hand the
estimates provided by Theorem \ref{lsthm2} is potentially stronger than \cite{F1}.

In this paper we also give a cohomological lower bound for $\Cat(X,\xi)$. 
It uses cup-products of cohomology classes 
of flat line bundles with monodromy described by complex numbers, which are not Dirichlet units.
Using this theorem one may get useful estimates for $\Cat(X,\xi)$ in specific examples.

In a previous paper \cite{F1} we introduced an invariant $\cat(X,\xi)$ which is similar to $\Cat(X,\xi)$.
We show that $\Cat(X,\xi)\geq \cat(X,\xi)$ always and  in some examples $\Cat(X,\xi) \not= \cat(X,\xi)$.

I would like to thank Shmuel Weinberger for useful discussions and help.

\section{A review of paper \cite{F1}}\label{review}

In this section we will briefly review a recent paper \cite{F1}. 
The results of the present paper represent a development of that earlier work \cite{F1}, and this review may help the reader 
to gain a comprehensive impression and 
make some comparisons. Besides, we will give in this section a new statement (Theorem \ref{lsmain5}), which was actually proven in
\cite{F1}, but stated in a weaker form (see Theorems \ref{lsmain4} and \ref{lsmain2}).
Some notions introduced in this section will be used in the rest of the article.

\subsection{Closed 1-forms on topological spaces}\label{lsapp1}

It will be convenient for us to use the notion of a
continuous closed 1-form on a topological space. In this subsection we will briefly describe this notion, referring to 
\cite{F1} for more details.

Let $X$ be a topological space.
{\it A continuous closed 1-form $\omega$ on $X$} 
is defined as a collection $\{f_U\}_{U\in \mathcal U}$
of continuous real-valued functions $f_U: U\to \R$, where $\mathcal U=\{U\}$ is an open cover of $X$,
such that for any pair $U, V\in \mathcal U$ the difference
\[f_U|_{U\cap V} - f_V|_{U\cap V}: U\cap V \to \R\]
is a locally constant function.  Another such collection $\{g_V\}_{V\in \mathcal V}$ 
(where $\mathcal V$ is another open cover of $X$) defines {\it an equivalent} closed 1-form 
if for any $U\in \mathcal U$ and $V\in \mathcal V$ the function 
$f_U -g_V: U\cap V\to \R$ is locally constant.

Consider the following exact sequence of sheaves over $X$
\begin{eqnarray}
0\to \R_X \to C_X\to B_X\to 0.\label{lsexact}
\end{eqnarray}
Here $\R_X$ denotes the sheaf of locally constant functions, $C_X$ denotes the sheaf of real-valued 
continuous functions, and
$B_X$ denotes the sheaf of germs of continuous functions modulo locally constant.
More precisely, $B_X$ is the sheaf corresponding to the presheaf $U\mapsto C_X(U)/\R_X(U)$. 
Comparing this with our definition of a continuous closed 1-form, we find that {\it the space of global sections of $B_X$,
$ H^0(X;B_X),$
coincides with the space of continuous closed 1-forms on $X$}. 

As an example, consider an open cover $\mathcal U=\{X\}$ consisting of the whole space $X$. Any continuous 
function $f: X\to \R$ defines a closed 1-form on $X$, which is denoted by $df$ (the differential of $f$). 

One may integrate continuous closed 1-forms along continuous paths.
Let $\omega$ be a continuous closed 1-form  on $X$ given 
by a collection of continuous real-valued functions $\{f_U\}_{U\in \mathcal U}$
with respect to an open cover $\mathcal U$ of $X$. Let $\gamma: [0,1]\to X$ be a continuous path. The line integral
$\int_\gamma \omega$ is defined as follows. Find a subdivision $t_0=0<t_1<\dots < t_N=1$ of the interval
$[0,1]$ such that for any $i$ the image $\gamma[t_i,t_{i+1}]$ is contained in a single open set $U_i\in \mathcal U$.
 Then we define
\begin{eqnarray}
\int\limits_\gamma \omega \, =\, 
\sum\limits_{i=0}^{N-1} \, \, [f_{U_i}(\gamma(t_{i+1})) - f_{U_i}(\gamma(t_i))].\label{lsintegral}
\end{eqnarray}
The standard argument shows that the integral (\ref{lsintegral}) does not depend on the choice of the subdivision and the 
open cover $\mathcal U$. 

\begin{lemma} Let $\omega$ be a continuous closed 1-form on topological spaces $X$.
For any pair of continuous
paths $\gamma, \gamma': [0,1]\to X$ with 
common beginning $\gamma(0)=\gamma'(0)$ and common endpoints $\gamma(1)=\gamma'(1)$, it holds
that $\int_\gamma \omega \, =\,  \int_{\gamma'}\omega,$
provided that $\gamma$ and $\gamma'$ are homotopic relative to the boundary. 
\end{lemma}

Any closed 1-form defines the {\it homomorphism of periods}
\begin{eqnarray}
\pi_1(X, x_0) \to \R, \quad [\gamma]\mapsto \int_\gamma \omega \in \R\label{lsperiod}
\end{eqnarray}
given by the integration of 1-form $\omega$ along closed loops $\gamma: [0,1]\to X$ with $\gamma(0)=x_0=\gamma(1)$. 
The homomorphism of periods (\ref{lsperiod}) is a group homomorphism. 

\begin{lemma} Let $X$ be a path-connected topological space. 
A continuous closed 1-form $\omega$ on $X$ equals $df$ for a
continuous function $f: X\to \R$ if and only if $\omega$ defines a trivial homomorphism of periods (\ref{lsperiod}). $\Box$

\end{lemma}

Any continuous 
closed 1-form $\omega$ on a topological space $X$ defines a {\it singular  cohomology class} $[\omega]\in H^1(X;\R)$.
It is defined by the homomorphism of periods (\ref{lsperiod}) viewed as an element of 
$\Hom(H_1(X);\R)=H^1(X;\R)$. As follows from the above Lemma, {\it two continuous
closed 1-forms $\omega$ and $\omega'$ on 
$X$ have the same cohomology class $[\omega]=[\omega']$
if and only if their difference
$\omega - \omega'$ equals $df$, where $f:X\to \R$ is a continuous function. }

From (\ref{lsexact}), using that $C_X$ is a fine sheaf, we obtain an exact sequence
\[0\to H^0(X;\R)\to H^0(X;C_X) \stackrel d\to H^0(X;B_X) \stackrel {[\, \, ]}\to H^1(X;\R_X)\to 0.\]
Here $H^0(X;C_X)=C(X)$ is the set of all continuous functions on $X$, and the map $d$ acts by 
assigning to a continuous function $f: X\to \R$ the closed 1-form $df\in H^0(X;B_X)$ (its differential). 
The group $H^1(X;\R_X)$ is the ${\rm \check C}$ech cohomology $\check H^1(X;\R)$ and the map 
$[\, \, ]$ assigns to a closed 1-form $\omega$ its $ \rm\check C$ech cohomology class $[\omega]\in \check H^1(X;\R)$. 

The natural map $\check H^1(X;\R) \to H^1(X;\R)$ is an isomorphism assuming that $X$ is paracompact, Hausdorff, 
and homologically locally connected, cf. \cite{Sp}. 
Recall that a topological space $X$ is {\it homologically locally connected} 
if for every point $x\in X$ and a neighborhood $U$ 
of $x$ there exists a neighborhood $V$ of $x$ in $U$ such that $\tilde H_q(V)\to \tilde H_q(U)$ is trivial for all $q$.

\begin{lemma} Let $X$ be a paracompact, Hausdorff, homologically locally connected topological space. 
Then any singular cohomology class $\xi\in H^1(X;\R)$ 
may be realized by a continuous closed 1-form on $X$. 
\end{lemma}

\subsection{A generalization of the Lusternik -- Schnirelman category}

Let $X$ be a finite CW-complex and let $\xi\in H^1(X;\R)$ be a real cohomology class. 
Fix a continuous closed 1-form $\omega$ on $X$ representing the cohomology class $\xi$. 

\begin{definition}\label{defcatxi1}
{\it We will define ${\rm cat}(X,\xi)$ 
to be the least integer $k$ such that for any integer $N>0$ there exists an open cover
\begin{eqnarray}
X=F\cup F_1\cup\dots \cup F_k,\label{cover1}
\end{eqnarray}
such that:
\begin{enumerate}
\item[(a)] Each inclusion $F_j\to X$ is null-homotopic, where $j=1, \dots, k$.
\item[(b)] There exists a homotopy $h_t: F\to X$, where $t\in [0,1]$, such that $h_0$ is the inclusion $F\to X$
and for any point $x\in F$,
\begin{eqnarray}
\int\limits_{\gamma_x}\omega \, \le \, -N,\label{lessthan1}
\end{eqnarray}
where the curve $\gamma_x:[0,1]\to X$ is given by $\gamma_x(t) =h_t(x)$.
\end{enumerate}}
\end{definition}

Intuitively, condition (b) means that in the process of the homotopy $h_t$ 
every point of $F$ makes at least $N$ full twists (in the negative direction)
with respect to $\omega$.

The main properties of $\cat(X,\xi)$ are:

\begin{itemize}

\item  $\cat(X,\xi)$ does not depend on the choice of continuous closed 
1-form $\omega$ (which appears in Definition \ref{defcatxi})
but only on the cohomology class $\xi=[\omega]$; 
\item $\cat(X,\xi)$ is a homotopy invariant of the pair $(X,\xi)$;
\item $\cat(X,\xi)=\cat(X)$ (the classical Lusternik -- Schnirelman category) for $\xi=0$;
\item if $X$ is connected and $\xi\neq 0$ then $\cat(X,\xi)\leq \cat(X)-1$; 
\item $\cat(X,\xi)$ admits a cohomological lower bound which uses cup-products of cohomology of generic flat line bundles 
(see \cite{F1}, Theorem 6.1) and also using the higher Massey products (see \cite{F1}, Theorem 6.4);
\item  in general, $\cat(X,-\xi)\neq \cat(X,\xi)$, see \cite{F1}, Example 3.4. 

\end{itemize}

The invariant $\cat(X,\xi)$ allows us to estimate the number of zeros of closed 1-forms in a given cohomology class:

\begin{theorem}{\rm (Theorem 4.1 in \cite{F1})}
\label{lsmain4} Let $\omega$ be a smooth closed 1-form on a closed manifold $M$ and let 
$\xi=[\omega]\in H^1(M;\R)$ denote the cohomology class of $\omega$.
If $\omega$ admits a gradient-like vector field $v$ with no homoclinic cycles, then $\omega$ has at least 
$\cat(M,\xi)$ geometrically distinct zeros.
\end{theorem}

Let $v$ denote a vector field on manifold $M$.
A {\it homoclinic orbit} of $v$ is defined as an integral trajectory $\gamma(t)$, $\dot\gamma(t)=v(\gamma(t))$, $t\in \R$,
 such that both limits $\lim_{t\to + \infty}\gamma(t)$ and $\lim_{t\to - \infty}\gamma(t)$ exist and are equal. 
A {\it homoclinic cycle} of $v$ 
is defined as a sequence of integral trajectories $\gamma_1, \gamma_2, \dots, \gamma_n$ of $v$
such that 
$\lim_{t\to +\infty} \gamma_i(t) = \lim_{t\to -\infty} \gamma_{i+1}(t)$
for $i=1, \dots, n-1$ and 
$\lim_{t\to +\infty} \gamma_n(t) = \lim_{t\to -\infty} \gamma_{1}(t).$ 
\begin{figure}[h]
\begin{center}
\includegraphics[0,0][310,108]{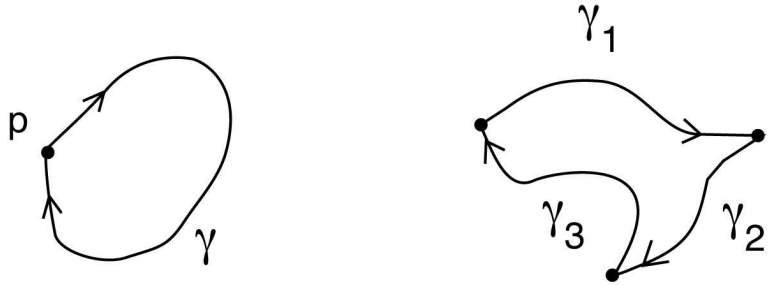}
\end{center}
\caption{Homoclinic orbit (left) and homoclinic cycle (right)}
\end{figure}

Here is a reformulation of the above Theorem:

\begin{theorem}{\rm (Theorem 4.2 in \cite{F1})}
\label{lsmain2} 
If the number of zeros of a smooth closed 1-form $\omega$ is less than $\cat(M,\xi)$, where
$\xi=[\omega]\in H^1(M;\R)$ denotes the cohomology class of $\omega$,
 then any
 gradient-like vector field for $\omega$ has a homoclinic cycle.
\end{theorem}

The proof of Theorem \ref{lsmain2}, given in \cite{F1}, actually proves the following slightly stronger statement:

\begin{theorem}\label{lsmain5}
Let $\omega$ be a closed 1-form on a closed manifold $M$ having less than $\cat(M,\xi)$ zeros, where
$\xi=[\omega]\in H^1(M;\R)$. Then there exists an integer $N>0$ such that any gradient-like vector field $v$ for $\omega$
has a homoclinic cycle $\gamma_1, \dots, \gamma_n$ with 
$$\sum_{i=1}^n\,  \int\limits_{\gamma_i}\omega \, \leq\,  N.$$
\end{theorem}

In other words, if closed 1-form 
$\omega$ has less than $\cat(M, \xi)$ zeros then any gradient like vector field $v$ for $\omega$ has 
a homoclinic cycle making at most $N$ full twists with respect to $\omega$.

\section{Definition of ${\rm Cat}(X,\xi)$}

In this section we will define invariant $\Cat(X,\xi)$ which is 
a modification of $\cat(X,\xi)$ described above. Later in this paper we will give a cohomological lower bound for it;
we will also use it to describe recurrence properties of dynamical systems and relations with the number of zeros.

Let $X$ be a finite CW-complex and let $\xi\in H^1(X;\R)$ be a real cohomology class. 

\begin{definition}\label{defcatxi}
{\it Fix a continuous closed 1-form $\omega$ on $X$ representing the cohomology class $\xi$, see \S \ref{review}.
Define ${\rm Cat}(X,\xi)$ 
to be the least integer $k$, such that there exists an open cover
\begin{eqnarray}
X=F\cup F_1\cup F_2\cup\dots \cup F_k,\label{cover}
\end{eqnarray}
with the following properties:
\begin{enumerate}
\item[(a)] Each inclusion $F_j\to X$, where $j=1, \dots, k$, is null-homotopic.
\item[(b)] There exists a homotopy $h_t: F\to X$, where $t\in (-\infty,+\infty)$, such that $h_0$ is the inclusion $F\to X$
and for any integer $n>0$ there exists $N>0$ such that for any $t>N$ and $x\in X$,
\begin{eqnarray}
\int\limits_x^{h_t(x)}\omega \, >\, n,\label{lsplus}
\end{eqnarray}
and for any $t<-N$ and $x\in X$,
\begin{eqnarray}
\int\limits_x^{h_t(x)}\omega \, <\,  -n.\label{lsminus}
\end{eqnarray}
\end{enumerate}}
\end{definition}

Intuitively, (\ref{lsplus}) and (\ref{lsminus}) can be expressed by saying that for any 
point $x\in F$,
\begin{eqnarray}
\lim_{t\to +\infty}\, \int\limits_{x}^{h_t(x)}\omega \, = \, +\infty,\quad
\lim_{t\to -\infty}\, \int\limits_{x}^{h_t(x)}\omega \, = \, -\infty
\label{lessthan}
\end{eqnarray}
and both limits are uniform in $x\in F$. 

 The meaning of the line integral $\int_\gamma \omega$, where $\gamma: [a, b]\to X$ is a continuous curve,
is explained in section \ref{review}.
In the left integral of (\ref{lessthan}) we integrate $\omega$ 
along the curve $\gamma_x: [0,t]\to X$, given by $\gamma_x(\tau) =h_{\tau}(x)$ (here we assume that $t>0$). 
We understand the right integral of (\ref{lessthan}) in a similar way: we integrate the closed 1-form $\omega$ along the curve
$[0, -t]\mapsto X$, where $\tau\mapsto h_{-\tau}(x)$ for $\tau\in [0, -t]$; here we assume that $t<0$.

Intuitively, condition (b) means that in the process of the homotopy $h_t$ 
every point of $F$ makes an infinite number of full twists in the positive (negative) direction
with respect to $\omega$ as $t$ tends to $+\infty$ (or, to $-\infty$, correspondingly). 

Observe that $\Cat(X,\xi)$ does not depend on the choice of the continuous closed 
1-form $\omega$ (which appears in Definition \ref{defcatxi})
but only on the cohomology class $\xi=[\omega]$. 
Indeed, if $\omega'$ is another continuous closed 1-form representing
$\xi$ then $\omega-\omega'=df$, where $f:X\to \R$ is a continuous function. Then for any continuous 
curve $\gamma:[a,b]\to X$, 
\[|\int_\gamma \omega - \int_\gamma \omega'| = |f(\gamma(b)) - f(\gamma(a))| \le C,\]
where the constant $C$ is independent of $\gamma$. Here we use the compactness of $X$. 
If one can construct an open cover (\ref{cover}) of $X$
such that condition (b) is satisfied, then the same is true for $\omega'$ replacing $\omega$. 

The following inequality holds
\begin{eqnarray}
\Cat(X,\xi) \le \cat(X).\label{usual}
\end{eqnarray}
 Here $\cat(X)$ denotes the
classical Lusternik -- Schnirelman category of $X$, i.e. the least integer $k$, such that
there is an open cover
$
X=F\cup F_1\cup\dots \cup F_k,
$
with
each inclusion $F_j\to X$ null-homotopic. Inequality (\ref{usual}) follows since we can always consider covers (\ref{cover}) with $F$ empty.
(\ref{usual}) can be improved, 
\begin{eqnarray}
\Cat(X,\xi) \le \cat(X)-1,\label{usual1}
\end{eqnarray}
assuming that $X$ is connected and $\xi\neq 0$. 
Under the above assumptions any open subset $F\subset X$, which is contractible to a point
in $X$, satisfies condition (b) of Definition \ref{defcatxi}. Hence, given a categorical open cover 
$X=G_1\cup \dots\cup G_r$ we may set
$F=G_1$ and $F_j =G_{j+1}$ for $j=1, \dots, r-1$, which gives a cover of $X$ with $r$ sets satisfying Definition \ref{defcatxi}.

Observe also that 
\begin{eqnarray}
\Cat(X, \xi) \, =\, \Cat(X, \lambda \xi), \quad \lambda\in \R,\label{lslambda}
\end{eqnarray}
as follows directly from the above definition. 

Let us compare the invariant $\Cat(X,\xi)$ with $\cat(X,\xi)$, see section \ref{review}. 
In Definition \ref{defcatxi} we make a stronger assumption on the subset $F$: we
require that it admits a deformation, rotating any point infinitely many times in {\it both directions}.
In the definition of $\cat(X,\xi)$ this condition was considerably weaker.
It is clear that  the following inequality 
\begin{eqnarray}\label{leq}
\Cat(X,\xi)\geq \cat(X,\xi).
\end{eqnarray}
holds.

In Example 3.4 of \cite{F1} we showed that it may happen that $\cat(X,-\xi)\not=\cat(X,\xi)$.
Since $\Cat(X,-\xi)=\Cat(X,\xi)$, see (\ref{lslambda}), we find that in this example, (\ref{leq}) cannot be an equality for both
classes $\xi$ and $-\xi$. We conclude that sometimes $\Cat(X,\xi) > \cat(X,\xi)$ and  hence {\it in general the invariants $\cat(X,\xi)$ and 
$\Cat(X,\xi)$ are distinct.}

Note that $\Cat(X,\xi)$ coincides with the Lusternik -- Schnirelman category $\cat(X)$ for $\xi=0$. Indeed, if $\xi=0$ we may take 
$\omega=0$ as a representing closed 1-form and then condition (b) of Definition \ref{defcatxi} may hold only if $F=\emptyset$.

{\bf Example.} Let $Y$ be a finite polyhedron.
Consider a bouquet $X=Y\vee S^1$ and assume that the class $\xi\in H^1(X,\R)$ satisfies
$\xi|_Y =0$ and $\xi|_{S^1}\neq 0$. Let us show that then
\[{\rm Cat}(X,\xi) \, =\, {\rm cat}(Y)-1.\]
Consider an open cover $X = F\cup F_1\cup\dots \cup F_k$ satisfying conditions (a) and (b) of Definition \ref{defcatxi}.
Let $F'$ denote $F\cap Y$ and $F_j'=F_j\cap Y$. 
We want to show that $F'$ is contractible in $Y$. 
This would imply that $Y=F'\cup F'_1\cup\dots \cup F_k'$ is a categorical cover of $Y$, and hence $\Cat(X,\xi)+1\le \cat(Y)$.

Let $h_t: F\to X$, $t\in \R$, be a homotopy
as in (b) of Definition \ref{defcatxi}.
 The infinite cyclic cover $\tilde X$ of $X$ corresponding to class $\xi$ is a union of a real line (the universal 
covering of the circle $S^1$) and an infinite number of disjoint copies of $Y$, each having a single common point with the line. 
The homotopy $h_t|_{F'}: F'\to X$ lifts to a homotopy $\bar h_t: F'\to \tilde X$, $t\in \R$, where $\bar h_0$ is the 
inclusion of $F'$
into a fixed copy of $Y\subset \tilde X$. Let $r: \tilde X\to Y$ be the retraction of $\tilde X$ onto $Y$. For $t=0$,
the map $r\circ \bar h_t : F'\to Y$ is the inclusion $F'\to Y$ and for $t>N$, $r\circ \bar h_t$ 
is the constant map into the base point of $Y$. Hence $r\circ \bar h_t$ is a deformation of $F'$ into a constant map.

The opposite inequality $\cat(X,\xi)+1\geq \cat(Y)$ is clear since for any categorical open cover $Y=G_0\cup \dots \cup G_r$
we may set $F=S^1\cup G_0$ and $F_j=G_j$, where $j=1, \dots, r$. The set $F$ satisfies (b) of Definition \ref{defcatxi}.

\section{Homotopy invariance of $\Cat(X,\xi)$}

\begin{lemma} Let $\phi: X_1\to X_2$ be a homotopy equivalence, $\xi_i\in H^1(X_i;\R)$, where $i=1, 2$,
and $\xi_1=\phi^\ast(\xi_2)\in H^1(X_1;\R)$.
Then 
\begin{eqnarray}
{\rm Cat}(X_1,\xi_1)={\rm Cat}(X_2, \xi_2).
\end{eqnarray}
\end{lemma}
\begin{proof} Let $\psi:X_2\to X_1$ be a homotopy inverse of $\phi$. 
Choose a closed 1-form $\omega_1$ on $X_1$ in the 
cohomology class $\xi_1$. 
Then $\omega_2=\psi^\ast \omega_1$ is a closed 1-form on $X_2$ lying in the cohomology class $\xi_2$. 

Fix a homotopy $r_t: X_1\to X_1$, where $t\in [0,1]$, such that $r_0=\id_{X_1}$ and $r_1= \psi\circ \phi$.
Compactness of $X_1$ implies that there is a constant $C>0$ such that $|\int_{\alpha_x}\omega_1|< C$
for any point $x\in X_1$, where $\alpha_x$ is the track of the point $x$ under homotopy $r_t$, i.e. 
$\alpha_x(t)=r_t(x)$.

Suppose that ${\rm Cat}(X_2, \xi_2)\leq k$. There is an open covering 
$X_2=F\cup F_1\cup\dots\cup F_k$, such that $F_1, \dots, F_k$ are contractible in $X_2$ and 
there exists a homotopy $h_t: F\to X_2$, where $t\in (-\infty,\infty)$, such that for any $x\in F$ the integral 
$$\int_{x}^{h_t(x)} \omega_2$$ 
(calculated along the curve
$\tau \mapsto h_{\tau}(x)$) 
tends to $\pm \infty$
as $t\to \pm \infty$ uniformly in $x\in F$.

Define 
\[G=\phi^{-1}(F), \quad G_j=\phi^{-1}(F_j), \quad j=1, \dots, k.\]
These sets form an open cover of $X_1=G\cup G_1\cup\dots \cup G_k$.
Let us show that the set $G\subset X_1$ satisfies condition (b) of Definition
\ref{defcatxi}. Define a homotopy $h'_t: G\to X_1$, where $t\in (-\infty,\infty)$, by
\begin{eqnarray}
h'_t(x) = \left\{
\begin{array}{ll}
\psi(h_{t+1}(\phi(x)))&\mbox{for}\quad t\leq -1,\\ 
r_{-t}(x), &\mbox{for}\quad -1\leq t\leq 0,\\ 
r_{t}(x), &\mbox{for}\quad 0\le t\le 1,\\ 
\psi(h_{t-1}(\phi(x))), &\mbox{for}\quad 1\leq t < \infty.
\end{array}
\right.
\end{eqnarray} 
Then $h'_0$ is the inclusion $G\to X_1$ and for any point $x\in G$ holds 
$$\lim_{t\to \pm\infty}\, \int_{x}^{h'_t(x)}\, \omega_1 \, =\,  \pm \infty,$$ 
holds uniformly in $x\in G$, 
where the integral is calculated along the track of $x$ under homotopy $h_t'$. 

The following diagram 
\begin{eqnarray*}
\begin{array}{ccc}
G_j & \stackrel {\subset}\to & X_1\\
\phi \downarrow && \uparrow \psi\\
F_j & \stackrel{\subset}\to & X_2
\end{array}
\end{eqnarray*}
is homotopy commutative and the horizontal map below is null-homotopic. This shows that the inclusion $G_j\to X_1$
is null-homotopic, where $j=1, \dots, k$. 

The above argument proves that $\Cat(X_1,\xi_1) \leq \Cat(X_2, \xi_2)$. 

The inverse inequality follows similarly. $\Box$.
\end{proof}

\section{Cohomological lower bound for $\Cat(X,\xi)$ and Dirichlet units}

In this section we obtain a lower bound for $\Cat(X,\xi)$ using cup-products of cohomology classes of local systems, see
Theorem \ref{lscupa} below.

\subsection{Moving homology classes}

Let $\tilde X\to X$ be an infinite cyclic covering of a connected finite polyhedron $X$. 
Fix a homeomorphism $\tau: \tilde X \to \tilde X$ generating the group of covering translations. 
$\tilde X$ has two ends $e=\pm \infty$, see \cite{Fr}. Let $K\subset \tilde X$ be a compact subset such that 
$\tilde X$ is the union of translates $\tau^iK$, where $i\in \Z$ runs over all integers.
Then for any $r\in \Z$ the set $\cup_{i\geq r}\tau^iK$ is a neighbourhood of an end of $\tilde X$ which we will denote by $+\infty$.
The sets of the form $\cup_{i\leq r}\tau^iK$ are neighbourhoods of the other end which we will denote $-\infty$.
Clearly, these conventions do not depend on the choice of $K$.

\begin{definition}{\it We say that a homology class $z\in H_q(\tilde X;\Z)$ is movable to an end $e$ (where $e=\pm \infty)$ of $\tilde X$ 
if for any neighbourhood $U\subset \tilde X$ of $e$ there exists a cycle in $U$ representing $z$.}
\end{definition}

Note that the homology $H_q(\tilde X;\C)$ is a finitely generated module over the ring $\C[\tau, \tau^{-1}]$ of Laurent polynomials with complex coefficients.
Here $\tau\in \C[\tau, \tau^{-1}]$ acts on homology classes as the induced map of $\tau:\tilde X\to \tilde X$.

In the next Theorem we use the following notation. For a finitely generated $\C[\tau, \tau^{-1}]$-module $M$
and a complex number $\lambda$ we denote by $M_\lambda$ the set of all $m\in M$ such that 
$(\tau-\lambda)^km =0$ for some $k$. The submodules $M_\lambda$ are all zero except for finitely many.
Note that $M_\lambda$ with $\lambda=0$ is always zero.
The  $\C[\tau, \tau^{-1}]$-torsion submodule of $M$ is the direct sum $\oplus_\lambda M_\lambda$, 
where $\lambda\in \C$ 
runs over all complex numbers.

\begin{theorem}\label{thm-} Given a homology class $z\in H_q(\tilde X;\Z)$, the following conditions are equivalent:

\noindent
(a) there exists a nonzero integer $k\in \Z$, such that $kz$ is movable to the end $-\infty$ of the
infinite cyclic covering $\tilde X$;

\noindent
(b) the image of $z$ under the homomorphism $H_q(\tilde X;\Z)\to H_q(\tilde X;\C)$ belongs to the submodule
\[\bigoplus_{\lambda} (H_q(\tilde X;\C))_\lambda\subset H_q(\tilde X;\C),\]
where $\lambda\in \C$ runs over all algebraic integers.
\end{theorem}
\begin{proof} Assume that an integer multiple $z'=kz$ of a 
homology class $z\in H_q(\tilde X;\Z)$, where, $k\neq 0$,
is movable to $-\infty$. Realize $z'$ by a cycle $c$ in $\tilde X$
and specify a compact subset $K\subset \tilde X$ containing $c$ and such that $\tilde X= \bigcup_{j\in \Z} \tau^j K$. 
Fix an integer $N>0$, such that $\tau^N K$ is disjoint from $K$. 

Using Mayer -- Vietories sequences, we may easily conclude that 
the homology class $z'$ may be realized by a cycle in $\tau^r K$ for any $r\leq -N$. 

Consider
\[V = \bigcap_{r\leq -N} \im  [H_q(\tau^r K;\Z) \to H_q(\tilde X;\Z)]\, \subset H_q(\tilde X;\Z).\]
It is a finitely generated Abelian subgroup of $H_q(\tilde X;\Z)$.
Moreover, $V$ 
is invariant under $\tau$. Hence, 
there exists a monic integral polynomial $p(\tau)\in \Z[\tau]$, such that $p(\tau)$ acts trivially on $V$. Since $z'$ belongs to $V$,
we obtain $p(\tau)z'=0$. All roots of $p(\tau)=0$ are algebraic integers. This shows that the image of $z$ in $H_q(\tilde X;\C)$
belongs to the submodule $\oplus (H_q(\tilde X;\C))_\lambda$, where $\lambda$ runs over algebraic integers.

Conversely, suppose that the image of $z$ in $H_q(\tilde X;\C)$ belongs to the direct 
sum $\oplus_\lambda (H_q(\tilde X;\C))_\lambda$, where $\lambda$ runs over all algebraic integers. Then there exists a monic
polynomial $p(\tau)\in \Z[\tau]$ and a nonzero integer $k$, such that $kp(\tau)z = 0$ in $H_q(\tilde X;\Z)$.
Consider the ring $R = \Z[\tau]/J$, where $J$ is the ideal generated by $p(\tau)$. It is a free Abelian group of rank
$r=\deg(p)$, and the powers $\tau^i$, where
$i=0, 1, \dots, r=\deg(p)-1$, form an additive basis of $R$.
Hence for any $N>0$ there are integers $b_j\in \Z$ such that 
\[\tau^N = \sum_{j=0}^r b_j \tau^j\quad \mbox{in}\quad R.\]
Therefore 
\[ \tau^N z' = \sum_{j=0}^{r} b_j\tau^j z'\quad \mbox{and}\quad 
z' = \sum_{j=0}^{r} b_j\tau^{j-N} z'\quad \mbox{in}\quad H_q(\tilde X;\Z),\]
where $z'$ denotes $kz$.

Assume that $K\subset \tilde X$ is a compact such that 
$\tilde X = \bigcup_{j\in \Z}\tau^j K$
and class $z'$ can be realized by a cycle $c$ in $K$. The formula above shows that for any integer
$N$ class $z'$ may be realized by a cycle
$$\sum_{j=0}^{r} b_j\tau^{j-N} c \quad \mbox{in}\quad \bigcup_{0\leq j \leq r}\tau^{j-N} K.$$ 
Hence $z'$ is movable to the end $-\infty$ of $\tilde X$. $\Box$
\end{proof}

The following Theorem gives a similar criterion of movability of an integral homology class to the positive end $+\infty$
of $\tilde X$. 

\begin{theorem}\label{thm+} {\rm (1)} A nonzero integer multiple $kz$ of a homology 
class $z\in H_q(\tilde X;\Z)$, $k\in \Z$,
is movable to the end $+\infty$ of the
infinite cyclic covering $\tilde X$ if and only if 
 the image of $z$ under the coefficient homomorphism $H_q(\tilde X;\Z)\to H_q(\tilde X;\C)$ belongs to the submodule
\[\bigoplus_{\lambda} (H_q(\tilde X;\C))_\lambda\subset H_q(\tilde X;\C),\]
where $\lambda\in \C$ runs over all complex numbers, such that $\lambda^{-1}$ is an algebraic integer.

{\rm (2)} A nonzero integer multiple $kz$ of a homology 
class $z\in H_q(\tilde X;\Z)$
is movable to both ends $+\infty$ and $-\infty$ of the
infinite cyclic covering $\tilde X$ if and only if 
 the image of $z$ under the homomorphism $H_q(\tilde X;\Z)\to H_q(\tilde X;\C)$ belongs to the submodule
\[\bigoplus_{\lambda} (H_q(\tilde X;\C))_\lambda\subset H_q(\tilde X;\C),\]
where $\lambda\in \C$ runs over all Dirichlet units.
\end{theorem}

A proof of statement (1) of Theorem \ref{thm+} is obtained by replacing $\tau$ by $\tau^{-1}$ in the proof of Theorem \ref{thm-}.
Statement (2) follows from Theorem \ref{thm-} and from statement (1).

\subsection{Lifting property}

We will use the following notation.
Given a finite CW-complex $X$ and an integral cohomology class $\xi\in H^1(X;\Z)$, for any nonzero complex number $a\in \C$ 
there is a local system 
over $X$ with fiber $\C$ such that the monodromy along any loop $\gamma$ in $X$ is 
multiplication by the number $a^{\langle \xi, \gamma\rangle}\in \C.$
Here $\langle \xi, \gamma\rangle\in \Z$ denotes the value of the class $\xi$ on the loop $\gamma$.
This local system will be denoted $a^\xi$. The cohomology of this local system $H^q(X;a^\xi)$ 
is a vector space of finite dimension.

\begin{theorem}\label{lifting}
 Let $X$ be a connected finite CW-complex, $\xi\in H^1(X;\Z)$ an integral nonzero cohomology class, 
and $p: \tilde X\to X$ the corresponding
infinite cyclic covering. Suppose that $F\subset X$ is a subset such that $\xi|_F=0$ and there exists a rotational homotopy 
$h_t:F\to X$, $t\in \R$, satisfying condition (b) of Definition \ref{defcatxi}. Then for any $a\in\C^\ast$, 
which is not a Dirichlet unit, 
the homomorphism 
\begin{eqnarray}
H^i(X,F;a^\xi) \to H^i(X;a^\xi)\label{homomo}
\end{eqnarray}
(we use the singular cohomology) is an epimorphism. \end{theorem}

\begin{proof} Our statement is equivalent to the claim that the restriction homomorphism $H^i(X;a^\xi)\to H^i(F;a^\xi)$ vanishes.
The latter is equivalent to the claim that the dual map of homology
groups $H_i(F;b^{\xi})\to H_i(X;b^\xi)$, vanishes, where $b=a^{-1}$.

Let $p: \tilde X\to X$ denote the infinite cyclic covering determined by class $\xi$. 
Symbol $\tilde F$ will denote $p^{-1}(F)\subset \tilde X$.
Furthermore, denote by $\tau: \tilde X\to \tilde X$
a generator of the covering translations group determined by the following condition: Let $\sigma:[0,1]\to  \tilde X$ be a 
path such that $\sigma(1) =\tau(\sigma(0))$; then $p\circ \sigma$ is a closed loop in $X$ and we require that
$\langle \xi, [p\circ \sigma]\rangle =+1$, where $[p\circ \sigma]\in H_1(X)$ is the corresponding homology class.

Consider the following commutative diagram:

$$
\begin{array}{ccccccccc}
H_i(\tilde F)& \stackrel{\tau-b}\to & H_i(\tilde F)&\to & H_i(F;b^\xi)&\stackrel \alpha\to &H_{i-1}(\tilde F)& \stackrel{\tau-b}\to  &H_{i-1}(\tilde F)\\
\downarrow && \downarrow j_\ast && \downarrow \beta && \downarrow && \downarrow \\
H_i(\tilde X) &\stackrel{\tau-b}\to & H_i(\tilde X)&\stackrel{\gamma}\to & H_i(X;b^\xi)&\to &H_{i-1}(\tilde X)& \stackrel{\tau-b}\to  &H_{i-1}(\tilde X)\end{array}
$$
The rows are well-known exact sequences and the vertical maps are induced by the inclusion $F\to X$. All homology 
groups, except for the middle column, are with coefficients in $\C$. 

Since we assume that $\xi|_F=0$, the set $\tilde F$ is homeomorphic to a disjoint union of infinitely many copies of $F$.
Hence $H_{i-1}(\tilde F)$ is a free $\C[\tau,\tau^{-1}]$-module, and homomorphism $\alpha$ vanishes.

Let $V=\oplus_\lambda (H_i(\tilde X))_\lambda$, where $\lambda$ runs over all Dirichlet units. The image of the homomorphism 
$H_i(\tilde F;\Z) \to H_i(\tilde X;\Z)$ consists of homology classes which are movable to both ends $\pm\infty$ of $\tilde X$. Hence by statement (2) of 
Theorem \ref{thm+}, 
the image $j_\ast$ in the diagram above belongs to $V$. 
Note that $V$ is a finite dimensional vector space invariant under operators $\tau$ and $\tau^{-1}$. 
Since we assume that $a$ is not a Dirichlet unit 
it follows that $\tau-b: V\to V$ is an isomorphism. This shows that $\im (j_\ast) \subset \im (\tau-b)$ in the diagram, and therefore
$\gamma\circ j_\ast=0$ and  $\beta=0$.
$\Box$
\end{proof}

\begin{remark}\label{correct1}{\rm
 In my paper \cite{F},
purely algebraic Proposition 3 is incorrect, as stated. It turns out that Corollary 2.6 of \cite{F} requires stronger assumptions. 
Theorem \ref{lifting} above is the corrected version of Corollary 2.6 of \cite{F}. In Theorem \ref{lifting} we impose much stronger assumptions on
the set $F$: we require that it admits a homotopy in $X$, infinitely rotating $F$ in both directions, see property (b)
of Definition \ref{defcatxi}.
Theorems 1 and 2 of \cite{F} should be replaced by Theorem \ref{lsthm2} (see below) combined with the estimate from below
on the number $\Cat(X,\xi)$, given by Theorem \ref{lscupa}.
}
\end{remark}

\subsection{A lower bound for $\Cat(X,\xi)$}

Let $X$ be a finite CW-complex and $\xi\in H^1(X;\Z)$. 
For any nonzero $a\in \C$ 
we have a local system $a^\xi$ over $X$ (see above).
Note that for $a=1$ the local system $a^\xi$ is the constant local system 
$\C$. If $a, b\in \C^\ast$ are two nonzero numbers, there is an isomorphism of local systems 
$a^{\xi}\otimes b^{\xi} \simeq (ab)^{\xi}.$
Hence we have a well-defined cup-product pairing
\begin{eqnarray*}
\cup: \quad H^q(X;a^\xi)\otimes H^{q'}(X;b^\xi) \to H^{q+q'}(X; (ab)^\xi).
\end{eqnarray*}

\begin{theorem}\label{lscupa} Assume that there exist cohomology classes 
$$u\in H^{q}(X;a^\xi),\quad v\in H^{q'}(X;b^\xi),\quad w_j\in H^{d_j}(X;\C),$$
where $j=1, \dots, r,$
such that:

{\rm (i)} $d_1>0, \dots, d_r>0$; 

{\rm (ii)} the cup-product
\begin{eqnarray}
0\not=\, \, u\cup v\cup w_1\cup w_2\cup\dots \cup w_r\in H^d(X;(ab)^\xi)\label{lscupb}
\end{eqnarray}
is nontrivial, where $d =q+q' +d_1 +\dots +d_r$;

{\rm (iii)} neither $a$ nor $b$ is a Dirichlet unit.
Then $\Cat(X,\xi)>r$.
\end{theorem}

Recall that a nonzero complex number $a\in \C$ is a Dirichlet unit if both $a$ and $a^{-1}$ are algebraic integers.
\begin{proof} 
We will assume below that the class $\xi$ is nonzero and indivisible.
For $\xi=0$ Theorem \ref{lscupa} is true and well known.

Suppose that $\Cat(X,\xi)\le r$. We want to show that then 
any cup-product (\ref{lscupb}) satisfying conditions
(i), (ii), (iii) of Theorem \ref{lscupa} vanishes. 
$\Cat(X,\xi)\le r$ means that there is an open cover 
$X=F\cup F_1\cup \dots \cup F_r$ satisfying properties 
(a) and (b) of Definition \ref{defcatxi}.
Each cohomology class $w_j\in H^{d_j}(X;\C)$ 
may be lifted to a relative cohomology class $\tilde w_j\in H^{d_j}(X,F_j;\C)$, where $j=1, \dots, r$.
This follows from the exact sequence
\[H^{d_j}(X,F_j;\C) \to H^{d_j}(X;\C) \to H^{d_j}(F_j;\C),\]
since $d_j>0$ and the second map is zero (since, the inclusion $F_j\to X$ is null-homotopic). 
We may represent the remaining open 
set $F$ as the union $F=F'\cup F^{\prime\prime}$ of two open subsets $F', F^{\prime\prime}\subset X$
such that $\xi|_{F'}=0$ and $\xi|_{F^{\prime\prime}}=0$; here we use the assumption that the class $\xi$ is integral, i.e.
$\xi\in H^1(X;\Z)$.
By Theorem \ref{lifting}
we may lift the cohomology class $u\in H^q(X;a^\xi)$ to a relative cohomology class
$\tilde u\in H^q(X, F'; a^\xi)$. 
Similarly, by Theorem \ref{lifting}
we may lift the class $v\in H^{q'}(X;b^\xi)$ to a relative cohomology class
$\tilde v\in H^{q'}(X, F^{\prime\prime}; b^\xi)$. 
Therefore the product (\ref{lscupb}) is obtained from a product
\[\tilde u\cup \tilde v \cup \tilde w_1\cup \dots \cup \tilde w_r\, \in H^\ast(X,X;(ab)^\xi) =0\]
(lying in the trivial group) by restricting to $X$. 
Hence any cup-product (\ref{lscupb}) must vanish. $\Box$

\end{proof}

\begin{remark} 
{\rm Theorem \ref{lscupa} may be compared with Theorem 6.1 from \cite{F1}. Because of inequality (\ref{usual}),
both Theorem \ref{lscupa} and Theorem 6.1 from \cite{F1} may be applied to give a lower bound for $\Cat(X,\xi)$.
There exist examples showing that either of these two theorems may give
better estimates. Apart from this, Theorem \ref{lscupa} has an obvious advantage since it explicitly
describes the set of forbidden monodromies (the set of Dirichlet units). In Theorem 6.1 from \cite{F1} the set of forbidden 
monodromies, ${\rm Supp}(X,\xi)$, depends on the pair $(X,\xi)$ and in specific examples, requires calculation.}
\end{remark}

\section{Dynamics of flows with Lyapunov 1-forms}\label{final}

Let $M$ be a closed smooth Riemannian manifold and let $v$ be a smooth vector field on $M$.
We will assume that $v$ has at most finitely many zeros.
The induced metric on $M$ is denoted $d(\cdot, \cdot)$.

\begin{definition} A subset $U\subset M$ is called $v$-convex if for any integral trajectory $x(t)$ of $v$
the set 
$\{t; \, x(t)\in U\}\subset  \R$ is convex.
\end{definition}

\noindent
Geometrically, this means that any integral trajectory of $v$, once leaving $U$, does not return to $U$ again.

Is it true that any neighbourhood of a zero $p\in M$, $v_p=0$,
contains a $v$-convex neighbourhood of $p$?
In other words, do there exist arbitrarily small $v$-convex neighbourhoods of $p$?

Here is one geometric reason why such a convex neighbourhood may not exist.
Suppose that there exists an integral trajectory $x(t)$, $\dot x=v(x)$, which is homoclinic to the zero $p$, 
i.e. the limits of $x(t)$, when $t$ tends to $+ \infty$ and to $-\infty$, both exist and coincide with $p$. 
Then clearly $p$ does not admit a small $v$-convex neighbourhood. 

In general, absence of an arbitrarily small $v$-convex neighbourhood of a zero $p\in M$, $v_p=0$, 
means that arbitrarily close to $p$ passes a 
trajectory which first moves far away from $p$ but eventually returns and passes very close to $p$.

\begin{definition}{\rm (see \cite{Sh}, Chapter 3)} A point $m\in M$ is called chain-recurrent with respect to the flow
$t\mapsto m\cdot t$ determined by the field $v$ if for any $\delta>0$, $T>0$ there exist points $y_0, y_1, \dots, y_k\in M$
and numbers $s_0, s_1, \dots, s_{k}$ such that 
$$ s_j>T, \quad d(m, y_0)<\delta, \quad d(y_k\cdot s_k, m)<\delta,$$
and 
$$d(y_j\cdot s_j, y_{j+1})<\delta, \quad\mbox{for}\quad  j=0, 1, \dots, k-1.$$
\end{definition}

Let $R_v$ denote the set of all chain recurrent points of the flow of $v$.

Note that the zeros of $v$ always belong to $R_v$. For us it will be important to know whether the zeros are isolated points
of the chain recurrent set $R_v$.

The following definition yields a generalization of the classical notion of a Lyapunov function (see \cite{Sh}, chapter 3);
usually our condition (b) does not appear in the case of functions.

\begin{definition} A Lyapunov 1-form for $v$ is a smooth closed 1-form $\omega$, $d\omega=0$, 
such that:

(a) on the complement of the set of zeros of $v$, $\omega(v) >0$;

(b) for any zero $p\in M$ of $v$ there exists a neighbourhood $U$ of $p$ such that $\omega|_U =df$ and 
the field $v|_U$ coincides with the gradient of the smooth function $f:U\to \R$ with respect to a Riemannian metric.
\end{definition}

\begin{lemma} \label{lemmarec}
Assume that $v$ admits a Lyapunov closed 1-form and the zeros of $v$ are isolated points of the chain-recurrent set
$R_v$. Then the zeros of $v$ admit arbitrarily small $v$-convex neighbourhoods.
\end{lemma}
\begin{proof} Let $p$ be one of the zeros of $v$, i.e. $v_p=0$. Assume that $p$ does not admit an arbitrarily small 
$v$-convex neighbourhood. 

Consider the universal cover $\tilde M$ of $M$ and the lift $\tilde v$ of $v$ into $\tilde M$. 
Note that $\tilde v$ is a gradient-like vector field on $\tilde M$.
Choose 
a pre-image $\tilde p\in \tilde M$ of $p$ and apply to it Lemma B.1 from \cite{F1}. 
We obtain a small $\tilde v$-convex neighbourhood $\tilde U$ of $\tilde p$. Denote by $\tilde A\subset \partial \tilde U$
the set of points $q\in \partial \tilde U$ such that the flow $q\cdot t$ tends to $\tilde p$ as $t\to +\infty$. Similarly, denote 
by $\tilde B\subset \partial \tilde U$
the set of points $q\in \partial \tilde U$ such that the flow $q\cdot t$ tends to $\tilde p$ as $t\to -\infty$. 
Let $U, A, B\subset M$ denote the images of the sets $\tilde U, \tilde A, \tilde B\subset \tilde M$ under the projection
$\tilde M\to M$.
Let $N$ denote the field of interior normals to $\partial U$. Denote
$$\partial_-U =\{q\in \partial U; (N_q, v_q)\leq 0\},\quad \partial_+U =\{q\in \partial U; (N_q, v_q)\geq 0\}.$$
Then $A\subset \partial_+U$, $B\subset \partial_-U$.

Assume that one cannot shrink the neighbourhood $U$ of $p$ to make it $v$-convex. Then there exists a sequence of points
$x_n\in \partial_-U$ and a sequence of numbers $t_n>0$ such that $x_n\cdot t_n\in \partial_+U$, 
$x_n$ converges to a point $b\in B$ and $x_n\cdot t_n$ converges to a point $a\in A$. Without loss of generality we may 
assume that the sequence $t_n$ has a finite or an infinite limit. If $t_n\to \tau\in \R$ then $b\cdot \tau =a$ and we see that there exists
an orbit homoclinic to $p$; any point of this orbit is chain recurrent, and so $p$ is not an isolated point of $R_v$.

Assume now that $\lim t_n =+\infty$.
We claim that in this case any point $q=b\cdot \sigma$ with $\sigma<0$ is chain recurrent. 
Suppose that $\delta>0$ and $T>0$ are given. Find $n$ so large that 
$$ d(x_n\cdot t_n, a)<\delta, \quad t_n>2T, \quad\mbox{and}\quad d(x_n\cdot T, b\cdot T)<\delta. $$
A recurrent $(\delta, T)$-chain starting and ending at $q$ can be defined as follows:
\begin{eqnarray*}
\begin{array}{lll}
y_0=q,& s_0=T-\sigma,&\\
y_1=x_n\cdot T,&s_1= t_n-T,\\
y_2=a, & s_2, &\mbox{$s_2$ is so large that $d(a\cdot s_2,p)<\delta/2$,}\\
y_3=b\cdot(-\alpha),& s_3=\sigma+\alpha,&\mbox{$\alpha$ is so large that $d(b\cdot (-\alpha),p)<\delta/2$}.
\end{array}
\end{eqnarray*}
This shows that arbitrarily close to $p$ we may find chain recurrent points distinct from $p$. This contradicts our assumption
that the zeros of $v$ are isolated points in the set of chain recurrent points.
 $\Box$
\end{proof}

In this section we will give a condition (in terms of the number of zeros of $v$ and the global properties of $v$)
which implies that the zeros of $v$ are not isolated in the set of chain recurrent points.

For $\epsilon >0$ we will denote by 
$U_\epsilon$ the set of all points $m\in M$ such that the integral trajectory
$x(t)=m\cdot t$ of $v$ satisfies the following: for all $t\in \R$ and all zeros $p$ ($v_{p}=0$) of $v$, $d(x(t),p)>\epsilon$. In other words, 
a point $m\in M$ belongs to $U_\epsilon$ if the trajectory through $m$ does not approach the zeros of $v$
by distance $\epsilon$.

Our additional assumption about vector field $v$ will be:

{\bf ($\ast$)} {\it for any $\delta>0$ there exists $ 0<\epsilon <\delta $ such that $U_\epsilon$ 
is a Euclidean Neighbourhood Retract (ENR).}

It is clear that property ($\ast$) of $v$ does not depend on the choice of Riemannian metric.

Recall that a topological space $X$ is an ENR if $X$ is homeomorphic to a subset $Y\subset \R^n$ such that $Y$ is a retract
of some open neighbourhood $\R^n\supset V\supset Y$, see \cite{D}, chapter 4.

\begin{theorem}\label{lsthm1} Let $v$ be a smooth vector field on a closed smooth manifold $M$ and 
let $\omega$ be a Lyapunov
1-form for $v$. Suppose that $v$ satisfies ($\ast$) and 
has less than $\Cat(M, \xi)$ zeros, where $\xi=[\omega]\in H^1(M;\R)$ denotes the 
cohomology class of $\omega$. Then at least one of the zeros of $v$ is not isolated in the chain recurrent set $R_v$ of $v$.
\end {theorem}

Here is an equivalent statement:

\begin{theorem}\label{lsthm2} Let $v$ be a smooth vector field on a closed smooth manifold $M$.
Assume that $v$ satisfies condition ($\ast$), 
admits a Lyapunov 1-form, and the zeros of $v$ are isolated points of the chain recurrent set $R_v$.
Then the number of zeros of $v$ is at least $\Cat(M,\xi)$, where $\xi\in H^1(M;\R)$ is any cohomology class 
which can be realized by a smooth
Lyapunov 1-form for $v$.
\end {theorem} 

\begin{proof} We will only prove Theorem \ref{lsthm2}, which is obviously equivalent to Theorem \ref{lsthm1}.
We will assume that $v$ has at most finitely many zeros, since otherwise the claim of the Theorem is clear.

Let $p_1, \dots, p_k\in M$ denote all the zeros of $v$. 
We will construct an open cover 
$$F\cup F_1\cup F_2\cup\dots \cup F_k=M,$$
such that for $j=1, \dots, k$ each inclusion $F_j\to M$ is null-homotopic, and the set $F$ 
admits a homotopy $h_t: V\to M$, where $t\in (-\infty, \infty)$,  infinitely rotating $F$ in $M$ with respect to $\omega$, as 
explained in Definition \ref{defcatxi}. This would imply that 
$$\Cat(M, \xi)\leq k,$$ 
confirming our claim. 

Denote by $M\times \R\to M$, $m\mapsto m\cdot t$, the flow generated by the vector field $v$.

For each zero $p_i$ of $v$ fix a small $v$-convex open neighbourhood $V_i$ around $p_i$ 
(it exists by Lemma \ref{lemmarec}).
We will assume that each $V_i$ is so small that the conclusion of Lemma B.1 from \cite{F1} applies to a lift of
$\bar V_i$ and of $v$ into the universal cover $\tilde M$. Besides, we will assume that $\bar V_i$ is contractible and 
$V_i\cap V_j =\emptyset$ for $i\not=j$. Fix a Riemannian metric on $M$ and let $\delta>0$ be such that
the closed $\delta$-ball around each zero $p_j$ is contained in $V_j$. 

Let $\epsilon>0$ be such that the set $U_\epsilon$ (see ($\ast$) above) is an ENR, $\epsilon <\delta$.
The manifold $M$ is clearly an ENR as well.
Using Corollary 8.7 from chapter 4 of \cite{D}, we find that there is an open neighbourhood $U\subset M$ of $U_\epsilon$
and a retraction $r: U\to U_\epsilon$, such that the inclusion $j: U\to M$ is homotopic to $i\circ r$, 
where $i: U_\epsilon\to M$ denotes the inclusion. Let $g_t: U\to M$, $t\in [0,1]$, be a homotopy with $g_0=j$ and $g_1=i\circ r$. 

Find an open set $F\supset V_\epsilon$ such that $\bar F\subset U$. 
Using the compactness of $\bar F$ we find that there exists a constant $C>0$ such that for all $x\in F$,
\begin{eqnarray}
|\int\limits_x^{r(x)}\omega |\, <\,  C,\label{boundls}
\end{eqnarray}
where the integral is calculated along the curve $t\mapsto g_t(x)$, $t\in [0,1]$.

Now we may define a homotopy
$h_t: F\to M$, where $t\in (-\infty, \infty)$, as follows:
\begin{eqnarray}
h_t(x) = \left\{
\begin{array}{lcc}
g_t(x), & \mbox{for} &\, 0\leq t\leq 1,\\
r(x)\cdot(t-1),& \mbox{for}&\,  t\geq 1,\\
g_{-t}(x), &\mbox{for} &\, -1\leq t\leq 0,\\
r(x)\cdot (t+1), & \mbox{for}&\,  t\leq -1.
\end{array}
\right.
\end{eqnarray}
This homotopy satisfies 
\begin{eqnarray}
\lim_{t\to \pm \infty} \int\limits_x^{h_t(x)} \omega \, = \, \pm \infty, \quad x\in F,\label{lsuniform}
\end{eqnarray}
and the limits are uniform with respect to $x\in F$. This easily follows from (\ref{boundls}) and from 
the following observation: There exists a constant $c>0$ such that $\omega(v)>c$ on the complement of the union of
$\epsilon$-balls around the zeros $p_j$. Note also that any trajectory starting at a point of $F$ does not
approach the zeros closer than $\epsilon$. Hence we find,
$$
\int\limits_x^{h_t(x)}\omega \, > \, ct -C, \quad t>0,\quad x\in F, 
$$
and 
$$
\int\limits_x^{h_t(x)}\omega \, < \, -ct +C, \quad t<0,\quad x\in F,
$$
which explain uniform convergence in (\ref{lsuniform}).

Let $F_i$ denote the set of all points $m\in M$ such that $m\cdot t$ belongs to $V_i$ for some $t$. 
Clearly, $M=F\cup F_1\cup \dots \cup F_k$, since any trajectory $x(t)=m\cdot t$ either enters a set $V_i$ (and then $m\in F_i$), 
or does not approach the zeros of the field $v$ by a distance less than $\delta$ (and then $m$ belongs to $F$). 

Let us show that $F_i$ is contractible in $M$. For a point $m \in F_i$, the set $\{t\in \R; m\cdot t \in V_i\}$ coincides with an open interval
$(a_m, b_m)$, where $a_m<b_m$ (because of $v$-convexity of $V_i$). 
The functions $m\mapsto a_m$ and $m\mapsto b_m$ are continuous real-valued functions on $F_i$ (see Lemma 
B.1 of \cite{F1}). Let $\psi_i: F_i\to \R$ be defined by 
$$\psi_i(m) = \mu(a_m, b_m), \quad m\in F_i,$$
where $\mu: \{(x,y)\in \R^2; x<y\}\to \R$
is given by
\begin{eqnarray*}
\mu(x,y) =\left\{
\begin{array}{lll}
0,&\mbox{if}& x\leq 0\,\,  \mbox{and}\,\,  y\geq 0,\\
x,&\mbox{if}& x\geq 0,\\
y,&\mbox{if}& y\leq 0.
\end{array}
\right.
\end{eqnarray*}
$\mu$ is continuous and therefore $\psi_i$ is continuous.
Consider the homotopy
\begin{eqnarray*}
f^i_t:\, F_i\to M, \quad 
f^i_t(m)=m\cdot(t\psi_i(m)), \quad m\in F_i, \quad t\in [0,1].
\end{eqnarray*}
Then $f^i_0(m)=m$ and $f^i_1(m)$ belongs to $\bar V_i$. Since $\bar V_i$ is assumed to be contractible, we obtain that
$F_i$ is contractible in $M$.

This completes the proof. $\Box$
\end{proof}

\bibliographystyle{amsalpha}

\vskip 2cm 

Address: 

Michael Farber,

School of Mathematical Sciences, 

Tel Aviv University, Ramat Aviv 69978, Israel

e-mail:   \, \, mfarber@tau.ac.il

\end{document}